\documentclass[11pt]{article}
\usepackage{amssymb,latexsym}
\usepackage{epsfig}
\usepackage{eufrak}
\usepackage{amsmath}
\usepackage{mathrsfs}
\usepackage{color}

\newtheorem{theorem}{Theorem}
\newtheorem{lemma}{Lemma}

\newcommand{\be}{\begin{equation}}
\newcommand{\ee}{\end{equation}}
\newcommand{\bee}{\begin{eqnarray*}}
\newcommand{\eee}{\end{eqnarray*}}
\newcommand{\bel}{\begin{eqnarray}}
\newcommand{\eel}{\end{eqnarray}}
\newcommand{\bec}{\begin{cases}}
\newcommand{\eec}{\end{cases}}
\newcommand{\bem}{\begin{bmatrix}}
\newcommand{\eem}{\end{bmatrix}}

\newcommand{\la}{\label}
\newcommand{\li}{\left}
\newcommand{\ri}{\right}

\newcommand{\lc}{\lceil}
\newcommand{\rc}{\rceil}
\newcommand{\lf}{\lfloor}
\newcommand{\rf}{\rfloor}

\newcommand{\vep}{\varepsilon}

\newcommand{\de}{\delta}

\newcommand{\se}{\theta}
\newcommand{\Se}{\Theta}

\newcommand{\al}{\alpha}
\newcommand{\ba}{\beta}

\newcommand{\Om}{\Omega}

\newcommand{\f}{\frac}

\newcommand{\cd}{\cdots}

\newcommand{\qu}{\quad}
\newcommand{\qqu}{\qquad}

\newcommand{\mscr}{\mathscr}

\newcommand{\bb}{\mathbb}

\newcommand{\wh}{\widehat}
\newcommand{\wt}{\widetilde}
\newcommand{\mrm}{\mathrm}
\newcommand{\bs}{\boldsymbol}

\newcommand{\tx}{\text}

\newcommand{\iy}{\infty}

\newcommand{\bed}{\begin{description}}
\newcommand{\eed}{\end{description}}
\newcommand{\bei}{\begin{itemize}}
\newcommand{\eei}{\end{itemize}}
\newcommand{\ben}{\begin{enumerate}}
\newcommand{\een}{\end{enumerate}}
\newcommand{\bib}{\bibitem}
\newcommand{\beL}{\begin{lemma}}
\newcommand{\eeL}{\end{lemma}}
\newcommand{\beT}{\begin{theorem}}
\newcommand{\eeT}{\end{theorem}}
\newcommand{\sect}{\section}

\newcommand{\bpf}{\begin{pf}}
\newcommand{\epf}{\end{pf}}
\newcommand{\bsk}{\bigskip}

\newcommand{\pfbox}{\hfill\mbox{$\Box$}}

\newenvironment{pf}{\paragraph*{Proof{\rm.}}}{\pfbox\bigskip}

\begin{document}

\title{{\bf Exact Sample Size Methods for Estimating Parameters of Discrete Distributions}
\thanks{
Dr. Xinjia Chen had been previously working with Louisiana State University at Baton Rouge, LA 70803, USA,  and is now with Department of
Electrical Engineering, Southern University and A\&M College, Baton Rouge, LA 70813, USA; Email: chenxinjia@gmail.com.  Dr. Zhengjia Chen is
working with Department of Biostatistics and Bioinformatics, Emory University, Atlanta, GA 30322; Email: zchen38@emory.edu.}}

\author{Xinjia Chen and Zhengjia Chen}

\date{November 2012}

\maketitle

\begin{abstract}

In this paper, we develop an approach for the exact determination of the minimum sample size for estimating the parameter of an integer-valued
random variable, which is parameterized by its expectation.  Under some continuity and unimodal property assumptions, the exact computation is
accomplished by reducing infinite many evaluations of coverage probability to finite many evaluations.  Such a reduction is based on our
discovery that the minimum of coverage probability with respect to the parameter bounded in an interval is attained at a discrete set of finite
many values.

\end{abstract}

\sect{Introduction}

Let $X$ be an integer-valued discrete random variable defined in a probability space $(\Om, \mscr{F}, \Pr)$ such that the probability mass
function is parameterized by its expectation
\[
\bb{E} [ X ] = \se \in \Se \subseteq (0, \iy), \]
 where $\Se$ is the parameter space.   It is a frequent problem to construct an estimator
$\wh{\bs{\se}}_n$ for $\se$ based on $n$ identical and independent samples $X_1, \cd, X_n$ of $X$.  An unbiased estimate of $\se$ is
conventionally taken as
\[
\wh{\bs{\se}}_n = \f{ Y_n } { n },  \] where
\[
Y_n = \sum_{i=1}^n X_i.
\]
  A crucial question in the estimation is as follows:

\bsk

{\it Given the knowledge that $\se$ belongs to interval $[a, b]$, what is the minimum sample size $n$ that guarantees the difference between
$\wh{\bs{\se}}_n$ and $\se$ be bounded within some prescribed margin of error with a confidence level higher than a prescribed value?}

\bsk

In this paper, we shall address this question based the assumption that {\it for any interval $\mscr{I}$, the probability $\Pr \li \{ Y_n \in
\mscr{I} \mid \se \ri \}$ is a continuous unimodal function of $\se \in \Se$}, where the notation $\Pr \{ E \mid \se \}$ denotes the probability
of event $E$ which is associated with the parameter $\se \in \Se$. This notation will be used throughout the paper. Clearly, the assumption is
satisfied if $X$ is a Bernoulli or Poisson random variable.

In this paper, the notion of a unimodal function is as follows:  A function $f(x)$ is said to be a unimodal function of $x \in [u, v]$ if there
exists a number $x^* \in [u, v]$ such that for any $x_1, x_2, x_3, x_4$ with $u \leq x_1 \leq x_2 \leq x^* \leq x_3 \leq x_4 \leq v$,
\[
f(u) \leq f(x_1) \leq f(x_2) \leq f(x^*), \qqu f(x^*)  \geq f(x_3) \geq f(x_4) \geq f(v).
\]

The paper is organized as follows. In Section 2, the techniques for computing the minimum sample size is developed with the margin of error
taken as a bound of absolute error. In Section 3, we derive corresponding sample size method by using relative error bound as the margin of
error.  In Section 4, we develop techniques for computing minimum sample size with a mixed error criterion.  In Section 5, we consider the
sample size problem in the context of range-preserving estimator.  Section 6 is the conclusion. The proofs are given in Appendices.   This work
is an extension of the recent works \cite{ChenJP, ChenPos} and \cite{Gamrot}.

Throughout this paper, we shall use the following notations. The set of integers is denoted by $\bb{Z}$.  The ceiling function and floor
function are denoted respectively by $\lc . \rc$ and $\lf . \rf$ (i.e., $\lc x \rc$ represents the smallest integer no less than $x$; $\lf x
\rf$ represents the largest integer no greater than $x$).  For integers $k \leq l$, the probability $\Pr \li \{ k \leq Y_n \leq l \mid \se \ri
\}$ is denoted by $S(n, k, l, \se)$. The left limit as $\eta$ tends to $0$ is denoted as $\lim_{\eta \downarrow 0}$.  The other notations will
be made clear as we proceed.

\section{Control of Absolute Error}

Let $\vep \in (0,1)$ be the margin of absolute error and $\de \in
(0,1)$ be the confidence parameter.  In many applications, it is
desirable to find the minimum sample size $n$ such that \[
 \Pr \li
\{ |\wh{\bs{\se}}_n - \se| < \vep \mid \se \ri \} > 1 - \de \] for any $\se \in [a,b]$. Here $\Pr  \{ |\wh{\bs{\se}}_n - \se | < \vep \mid \se
\}$ is referred to as the coverage probability.  The interval $[a, b]$ is introduced to take into account the knowledge of $\se$.  The exact
determination of minimum sample size is readily tractable with modern computational power by taking advantage of the behavior of the coverage
probability characterized by Theorem \ref{thm_abs} as follows.

\bsk

\beT \la{thm_abs} Let $\vep > 0$ and $0 \leq a < b$ such that $[a, b] \subseteq \Se$.  Assume that for any interval $\mscr{I}$, the probability
$\Pr \li \{ Y_n \in  \mscr{I} \mid \se \ri \}$ is a continuous unimodal function of $\se \in \Se$. Then, the minimum of $\Pr \{ |
\wh{\bs{\se}}_n - \se | < \vep \mid \se \}$ with respect to $\se \in [a, b]$ is achieved at the finite set $\{a, b \} \cup
 \{ \f{\ell}{n} + \vep \in (a, b) : \ell \in \bb{Z} \} \cup
 \{ \f{\ell}{n} - \vep \in (a, b) : \ell \in \bb{Z} \} $,
 which has less than $2n(b-a) + 4$ elements.
\eeT

See Appendix A for a proof.  The application of Theorem \ref{thm_abs} in the computation of minimum sample size is obvious.  For a fixed sample
size $n$, since the minimum of coverage probability with $\se \in [a, b]$ is attained at a finite set, it can determined by a computer whether
the sample size $n$ is large enough to ensure $\Pr \{ |\wh{\bs{\se}}_n - \se| < \vep \mid \se \} > 1 - \de$ for any $\se \in [a, b]$. Starting
from $n = 2$, one can find the minimum sample size by gradually incrementing $n$ and checking whether $n$ is large enough.

\section{Control of Relative Error}

Let $\vep \in (0,1)$ be the margin of  relative error and $\de \in
(0,1)$ be the confidence parameter.  It is interesting to determine
the minimum sample size $n$ so that
\[
\Pr \li \{  \li | \f{ \wh{\bs{\se}}_n - \se } {\se} \ri | < \vep \mid \se \ri \}
> 1 - \de
\]
for any $\se \in [a,b]$.  As has been pointed out in Section 2, an essential machinery is to reduce infinite many evaluations of the coverage
probability $\Pr \{ | \wh{\bs{\se}}_n - \se | < \vep \se \mid \se \}$ to finite many evaluations. Such reduction can be accomplished by making
use of Theorem \ref{thm_rev} as follows.

 \beT \la{thm_rev} Let $0 < \vep < 1$ and $0 < a  < b$ such that $[a, b] \subseteq \Se$. Assume that for any interval $\mscr{I}$, the probability
$\Pr \li \{ Y_n \in  \mscr{I} \mid \se \ri \}$ is a continuous unimodal function of $\se \in \Se$.
 Then, the minimum of $\Pr \li \{  \f{ |\wh{\bs{\se}}_n - \se| } {\se} < \vep \mid \se \ri \}$ with respect
to $\se \in [a, b]$ is achieved at the finite set $\{a, b \} \cup
 \{ \f{\ell}{n(1 + \vep)} \in (a, b) : \ell \in \bb{Z} \} \cup
 \{ \f{\ell}{n (1 - \vep)} \in (a, b) : \ell \in \bb{Z} \}$,
 which has less than $2n(b-a) + 4$ elements.
\eeT

See Appendix B for a proof.

\section{Control of Absolute Error or Relative Error}

Let $\vep_a > 0$ and $\vep_r \in (0,1)$ be respectively the margins of absolute error and relative error.  Let $\de \in (0,1)$ be the confidence
parameter.  In many situations, it is desirable to find the smallest sample size $n$ such that \be \la{rule2} \Pr \li \{ |\wh{\bs{\se}}_n - \se|
< \vep_a \; \; \mrm{or} \;\;  \li | \f{ \wh{\bs{\se}}_n - \se } {\se} \ri | < \vep_r \mid \se \ri \}
> 1 - \de \ee for any $\se \in
[a,b]$.  To make it possible to compute exactly
 the minimum sample size associated with (\ref{rule2}), we have
Theorem \ref{thm_abs_rev} as follows.

\beT \la{thm_abs_rev} Let $\vep_a > 0, \; 0 < \vep_r < 1$ and $0 \leq a < \f{\vep_a}{\vep_r} < b$ such that $[a, b] \subseteq \Se$. Assume that
for any interval $\mscr{I}$, the probability $\Pr \li \{ Y_n \in  \mscr{I} \mid \se \ri \}$ is a continuous unimodal function of $\se \in \Se$.
Then, the minimum of $\Pr \li \{ |\wh{\bs{\se}}_n - \se| < \vep_a \; \; \mrm{or} \;\; \li | \f{ \wh{\bs{\se}}_n - \se } {\se} \ri | < \vep_r
\mid \se \ri \}$ with respect to $\se \in [a, b]$ is achieved at the finite set $\{a, b, \f{\vep_a}{\vep_r} \} \cup
 \{ \f{\ell}{n} + \vep_a \in (a, \f{\vep_a}{\vep_r}) : \ell \in \bb{Z} \} \cup
 \{ \f{\ell}{n} - \vep_a \in (\f{\vep_a}{\vep_r}, b) : \ell \in \bb{Z} \} \cup
 \{ \f{\ell}{n(1 + \vep_r)} \in (a, \f{\vep_a}{\vep_r}) : \ell \in \bb{Z} \} \cup
 \{ \f{\ell}{n (1 - \vep_r)} \in (\f{\vep_a}{\vep_r}, b) : \ell \in \bb{Z}
 \}$, which has less than $2n(b-a) + 7$ elements.
\eeT

\bsk

Theorem \ref{thm_abs_rev} can be shown by applying Theorem
\ref{thm_abs} and Theorem \ref{thm_rev} with the observation that
\[
\Pr \li \{ |\wh{\bs{\se}}_n - \se| < \vep_a \; \; \mrm{or} \;\; \li | \f{ \wh{\bs{\se}}_n - \se } {\se} \ri | < \vep_r \ri \} = \bec \Pr \li \{
|\wh{\bs{\se}}_n - \se| < \vep_a  \ri \} & \tx{for} \; \se \in \li [ a,
\f{\vep_a}{\vep_r} \ri ], \\
\Pr \li \{ \li | \f{ \wh{\bs{\se}}_n - \se } {\se} \ri | < \vep_r \ri \} & \tx{for} \; \se \in \li ( \f{\vep_a}{\vep_r}, b \ri ]. \eec
\]

\bsk

\sect{Error Control for Range-Preserving Estimator}

In many situations, it may be appropriate to use the range-preserving estimator $\wt{\bs{\se}}_n$ for $\se$, which is defined as
\[
\wt{\bs{\se}}_n = \bec \wh{\bs{\se}}_n & \tx{for} \; \wh{\bs{\se}}_n \in [a, b],\\
a & \tx{for} \; \wh{\bs{\se}}_n < a,\\
b & \tx{for} \; \wh{\bs{\se}}_n > b. \eec
\]
See \cite{Gamrot} and the references therein. Recently, Gamrot \cite{Gamrot} has established an exact sample size method for estimating a
binomial parameter using the range preserving estimator.  Inspired by the work of  Gamrot \cite{Gamrot}, we have developed an exact approach for
the general sample size problem of estimating the parameter of a discrete distribution under certain unimodal property and continuity
assumptions. In the sequel, let $(u, v)$ denote an open interval if $u < v$. In the case that $u \geq v$, $(u, v)$ is an empty set.

To determine the minimum sample size for controlling absolute error in the context of using the range-preserving estimator, we have the
following results.

\beT

\la{Chensabs}

Assume that  $\vep > 0$ and $0 < a  < b < \iy$ such that $[a, b] \subseteq \Se$. Assume that for any interval $\mscr{I}$, the probability $\Pr
\li \{ Y_n \in \mscr{I} \mid \se \ri \}$ is a continuous unimodal function of $\se \in \Se$. Then, the minimum of $\Pr \{ | \wt{\bs{\se}}_n -
\se | < \vep  \mid \se \}$ with respect to $\se \in [a, b ]$ is attained at the finite set $A \cap [a, b]$, where {\small \be \la{Aabs} A = \li
\{ a, b, a + \vep, b - \vep \ri \} \bigcup \li \{  \f{k}{n} - \vep \in \li ( a, b -  \vep \ri ): k \in \bb{Z} \ri \} \bigcup \li \{ \f{k}{n} +
\vep \in \li ( a + \vep, b \ri ): k \in \bb{Z} \ri \},   \ee} which has less than $2 n (b - a - \vep)  + 6$ elements.  \eeT

See Appendix \ref{Chensabs_app} for a proof of Theorem \ref{Chensabs}.

To determine the minimum sample size for controlling relative error in the context of using the range-preserving estimator, we have the
following results.

\beT

\la{Chens}

Assume that $\vep \in (0, 1)$ and $0 < a  < b < \iy$ such that $[a, b] \subseteq \Se$. Assume that for any interval $\mscr{I}$, the probability
$\Pr \li \{ Y_n \in \mscr{I} \mid \se \ri \}$ is a continuous unimodal function of $\se \in \Se$. Then,  the minimum of $\Pr \{ |
\wt{\bs{\se}}_n - \se | < \vep \; \se \mid \se \}$ with respect to $\se \in [a, b ]$ is attained at the finite set $A \cap [a, b]$, where
{\small \be \la{Arev} A =  \li \{ a, b, \f{a}{1 - \vep}, \f{b}{1 + \vep} \ri \} \bigcup \li \{  \f{k}{n (1 + \vep)} \in \li ( a, \f{b}{1 + \vep}
\ri ): k \in \bb{Z} \ri \} \bigcup \li \{ \f{k}{n (1 - \vep)} \in \li ( \f{a}{1 - \vep}, b \ri ): k \in \bb{Z} \ri \}, \ee} which has less than
$2 n (b - a) - n \vep ( a + b) + 6$ elements.  \eeT

See Appendix \ref{Chens_app} for a proof.  It should be noted that Theorem \ref{Chens} is a generalization of Theorems 4 and 5 of Gamrot
\cite{Gamrot}.

To determine the minimum sample size for controlling absolute and relative error in the context of using the range-preserving estimator, we have
the following results.

 \beT \la{thm_abs_rev_range} Let $\vep_a > 0, \; 0 < \vep_r < 1$ and $0 \leq a  < \f{\vep_a}{\vep_r} < b$ such that $[a, b] \subseteq \Se$.
 Assume that for any interval $\mscr{I}$, the probability $\Pr \li \{ Y_n \in  \mscr{I} \mid \se \ri \}$ is a
continuous unimodal function of $\se \in \Se$. Then, the minimum of $\Pr \li \{ |\wt{\bs{\se}}_n - \se| < \vep_a \; \; \mrm{or} \;\; \li | \f{
\wt{\bs{\se}}_n - \se } {\se} \ri | < \vep_r \mid \se \ri \}$ with respect to $\se \in [a, b]$ is achieved at the finite set $\li ( A \cap [a,
\f{\vep_a}{\vep_r} ] \ri ) \cup \li ( B  \cap [\f{\vep_a}{\vep_r}, b] \ri )$, where {\small \bee  A =   \li \{ a,  a + \vep_a,
\f{\vep_a}{\vep_r} - \vep_a, \f{\vep_a}{\vep_r} \ri \} \bigcup \li \{ \f{k}{n} - \vep_a \in \li ( a, \f{\vep_a}{\vep_r} -  \vep_a \ri ): k \in
\bb{Z} \ri \} \bigcup \li \{ \f{k}{n} + \vep_a \in \li ( a + \vep_a, \f{\vep_a}{\vep_r} \ri ): k \in \bb{Z} \ri \}, \qu  &  & \\
B =  \li \{ b, \f{\f{\vep_a}{\vep_r}}{1 - \vep_r}, \f{b}{1 + \vep_r} \ri \} \bigcup \li \{  \f{k}{n (1 + \vep_r)} \in \li ( \f{\vep_a}{\vep_r},
\f{b}{1 + \vep_r} \ri ): k \in \bb{Z} \ri \} \bigcup \li \{ \f{k}{n (1 - \vep_r)} \in \li ( \f{\f{\vep_a}{\vep_r}}{1 - \vep_r}, b \ri ): k \in
\bb{Z} \ri \}.  \qu &   & \eee} Moreover, $A \cup B$ has less than $2n(b-a - \vep_a) - n (\vep_a + b \vep_r) + 11$ elements. \eeT

\bsk

Theorem \ref{thm_abs_rev_range} can be shown by applying Theorem \ref{Chensabs} and Theorem \ref{Chens} with the observation that
\[
\Pr \li \{ |\wt{\bs{\se}}_n - \se| < \vep_a \; \; \mrm{or} \;\; \li | \f{ \wt{\bs{\se}}_n - \se } {\se} \ri | < \vep_r \ri \} = \bec \Pr \li \{
|\wt{\bs{\se}}_n - \se| < \vep_a  \ri \} & \tx{for} \; \se \in \li [ a,
\f{\vep_a}{\vep_r} \ri ], \\
\Pr \li \{ \li | \f{ \wt{\bs{\se}}_n - \se } {\se} \ri | < \vep_r \ri \} & \tx{for} \; \se \in \li ( \f{\vep_a}{\vep_r}, b \ri ]. \eec
\]

\section{Conclusion}

We have developed an exact method for the computation of minimum sample size for estimating the parameter of an integer-valued discrete random
variable, which only requires finite many evaluations of the coverage probability. Our sample size method permits rigorous control of
statistical error for estimating parameters of common distributions such as binomial and Poisson distributions.

\appendix

\sect{Proof of Theorem \ref{thm_abs}}

Define \[ C(\se)  =  \Pr \li \{ \li | \f{Y_n}{n} - \se \ri | < \vep \mid \se \ri \} = \Pr \li \{ g(\se) \leq Y_n \leq h(\se) \mid \se \ri \}
\] where
\be \la{ghabs}
 g(\se) = \lf n( \se - \vep) \rf + 1, \qqu h(\se) = \lc n( \se + \vep) \rc - 1.
\ee It should be noted that $C(\se), \; g(\se)$ and $h(\se)$ are actually multivariate functions of $\se, \; \vep$ and $n$.  For simplicity of
notations, we drop the arguments $n$ and $\vep$ throughout the proof of Theorem \ref{thm_abs}.

We need some preliminary results.

\beL \la{minus} Let $\se_\ell = \f{\ell}{n} - \vep$ where $\ell \in \bb{Z}$. Then, $h(\se) = h(\se_{\ell + 1}) = \ell$ for any $\se \in
(\se_\ell, \se_{\ell +1})$. \eeL

\bpf For $\se \in ( \se_\ell, \; \se_{\ell + 1})$, we have $0 < n \li (\se - \se_\ell \ri ) < 1$ and \bee h (\se)  & = &
\lc n( \se + \vep) \rc - 1\\
& = &  \lc n( \se_\ell + \vep + \se - \se_\ell ) \rc - 1\\
& = & \li \lc n \li ( \f{\ell}{n} - \vep  + \vep + \se - \se_\ell  \ri )  \ri \rc - 1\\
& = & \ell - 1 + \li \lc n \li (\se - \se_\ell \ri )  \ri \rc\\
& = & \ell\\
& = & \li \lc n \li ( \f{\ell + 1}{n} - \vep  + \vep \ri )  \ri \rc - 1 = h(\se_{\ell + 1}). \eee

\epf

\beL \la{plus} Let $\se_\ell = \f{\ell}{n} + \vep$ where $\ell \in \bb{Z}$. Then, $g(\se) = g(\se_{\ell})= \ell + 1$ for any $\se \in (\se_\ell,
\se_{\ell +1})$. \eeL

 \bpf
For $\se \in \li ( \se_\ell, \; \se_{\ell + 1} \ri )$,  we have $-1 < n \li (\se - \se_{\ell + 1} \ri ) < 0$ and \bee
g (\se) & = & \lf n( \se - \vep) \rf + 1\\
& = & \lf n( \se_{\ell + 1} - \vep + \se - \se_{\ell + 1} ) \rf + 1\\
& = &  \li \lf n \li ( \f{ \ell + 1 } { n  } + \vep - \vep \ri ) \ri \rf +
\lf n( \se - \se_{\ell + 1} ) \rf +  1\\
& = & \li \lf n \li ( \f{ \ell + 1 } { n  } + \vep  - \vep \ri ) \ri \rf - 1 +  1\\
& = & \ell + 1\\
& = & \li \lf n \li ( \f{ \ell } { n  } + \vep - \vep \ri ) \ri \rf +  1 =  g(\se_\ell). \eee

\epf

\beL \la{constant} Let $\al < \ba$ be two consecutive elements of
the ascending arrangement of all distinct elements of $\{a, b \}
\cup
 \{ \f{\ell}{n} + \vep \in (a, b) : \ell \in \bb{Z} \} \cup
 \{ \f{\ell}{n} - \vep \in (a, b) : \ell \in \bb{Z} \} $.
Then, both $g(\se)$ and $h (\se)$ are constants for any $\se \in (\al, \ba)$.
 \eeL

 \bpf
Since $\al$ and $\ba$ are two consecutive elements of the ascending arrangement of all distinct elements of the set, it must be true that there
is no integer $\ell$ such that {\small $\al < \f{\ell}{n} + \vep < \ba$} or {\small $\al < \f{\ell}{n} - \vep < \ba$}.  It follows that there
exist two integers $\ell$ and $\ell^\prime$ such that {\small $(\al, \ba) \subseteq \li ( \f{\ell}{n} + \vep, \f{\ell + 1}{n} + \vep \ri )$} and
{\small $(\al, \ba) \subseteq \li ( \f{\ell^\prime}{n} - \vep, \f{\ell^\prime + 1}{n} - \vep \ri )$}. Applying Lemma \ref{minus} and Lemma
\ref{plus}, we have {\small $g(\se) = g \li ( \f{\ell}{n} + \vep \ri )$} and {\small $h(\se) = h \li (\f{\ell^\prime + 1}{n} - \vep \ri )$} for
any $\se \in (\al, \ba)$.

 \epf

\beL \la{lem_lim}
 For any $\se \in (0,1)$, $\lim_{\eta \downarrow 0} C(\se + \eta) \geq C(\se)$
 and $\lim_{\eta \downarrow 0} C(\se - \eta) \geq C(\se)$.
\eeL

\bpf

Observing that $h(\se + \eta) \geq h(\se)$ for any $\eta > 0$ and that
 \bee g(\se + \eta) & =  & \lf n( \se + \eta - \vep)
\rf + 1  \\
& = &  \lf n( \se - \vep) \rf + 1 + \lf n( \se - \vep) -
\lf n( \se  - \vep) \rf + n \eta\\
 & = & \lf n( \se - \vep)
\rf + 1  = g(\se) \eee for $0 < \eta < \f{ 1 + \lf n( \se - \vep) \rf -  n( \se - \vep)} {n}$, we have \be \la{ineqa} S(n, g(\se + \eta), h (\se
+ \eta), \se + \eta ) \geq S(n, g(\se), h (\se), \se + \eta ) \ee for $0 < \eta < \f{ 1 + \lf n( \se - \vep) \rf - n( \se - \vep)} {n}$. Since
 \[ h(\se + \eta)  =  \lc n( \se + \eta + \vep)
\rc - 1 = \lc n( \se + \vep) \rc - 1 + \lc n( \se + \vep) - \lc n( \se + \vep) \rc + n \eta \rc, \] we have
\[
h(\se + \eta) = \bec \lc n( \se + \vep) \rc & \tx{for} \; n( \se + \vep) = \lc n( \se  + \vep) \rc \; \tx{and} \; 0 < \eta <
\f{1}{n},\\
\lc n( \se + \vep) \rc  - 1 & \tx{for} \; n( \se + \vep) \neq \lc n( \se + \vep) \rc\; \tx{and} \; 0 < \eta < \f{\lc n( \se + \vep) \rc - n( \se
+ \vep)}{n}. \eec
\]
It follows that both $g(\se + \eta)$ and $h(\se + \eta)$ are independent of $\eta$ if $\eta > 0$ is small enough.  Since $S(n, g, h, \se +
\eta)$ is continuous with respect to $\eta$ for fixed $g$ and $h$, we have that $\lim_{\eta \downarrow 0} S(n, g(\se + \eta), h (\se + \eta),
\se + \eta )$ exists.  As a result, \bee \lim_{\eta \downarrow 0} C(\se + \eta) & = & \lim_{\eta \downarrow 0} S(n,
g(\se + \eta), h (\se + \eta), \se + \eta )\\
& \geq & \lim_{\eta \downarrow 0} S(n, g(\se), h (\se), \se + \eta ) = S(n,  g(\se), h (\se), \se ) = C(\se), \eee where the inequality follows
from (\ref{ineqa}).

Observing that $g(\se - \eta) \leq g(\se)$ for any $\eta > 0$ and that
 \bee h(\se - \eta) & =  & \lc n( \se - \eta + \vep)
\rc - 1\\
& = & \lc n( \se + \vep) \rc - 1 + \lc n( \se + \vep) -
\lc n( \se  + \vep) \rc - n \eta \rc\\
 & = & \lc n( \se + \vep)
\rc - 1   = h(\se) \eee for $0 < \eta < \f{ 1 + n( \se + \vep) - \lc n( \se + \vep) \rc } {n}$, we have \be \la{ineqb} S(n, g(\se - \eta), h
(\se - \eta), \se - \eta ) \geq S(n, g(\se), h (\se), \se - \eta ) \ee for {\small $0 < \eta < \min \li \{ \se, \f{ 1 + n( \se + \vep) - \lc n(
\se + \vep) \rc } {n} \ri \}$}.  Since \bee g(\se - \eta) & = & \lf n( \se - \eta - \vep)
\rf + 1 \\
& = & \lf n( \se - \vep) \rf + 1 + \lf n( \se - \vep) - \lf n( \se  - \vep) \rf - n \eta \rf, \eee we have {\small \[ g(\se - \eta) = \bec \lf
n( \se - \vep) \rf  & \tx{for} \; n( \se - \vep) = \lf n( \se  - \vep) \rf \; \tx{and} \; 0 < \eta
< \f{1}{n},\\
\lf n( \se - \vep) \rf  + 1 & \tx{for} \; n( \se - \vep) \neq \lf n( \se - \vep) \rf \; \tx{and} \; 0 < \eta < \f{n( \se - \vep) - \lf n( \se -
\vep) \rf }{n}. \eec
\]}
It follows that both $g(\se - \eta)$ and $h(\se - \eta)$ are independent of $\eta$ if $\eta > 0$ is small enough.  Since $S(n, g, h, \se -
\eta)$ is continuous with respect to $\eta$ for fixed $g$ and $h$, we have that $\lim_{\eta \downarrow 0} S(n, g(\se - \eta), h (\se - \eta),
\se - \eta )$ exists. Hence, \bee \lim_{\eta \downarrow 0} C(\se - \eta) & =
& \lim_{\eta \downarrow 0} S(n, g(\se - \eta), h (\se - \eta), \se -\eta )\\
& \geq & \lim_{\eta \downarrow 0} S(n, g(\se), h (\se), \se -\eta ) =  S(n, g(\se), h (\se), \se ) = C(\se), \eee where the inequality follows
from (\ref{ineqb}).

\epf

\beL \la{inbetween} Let $\al < \ba$ be two consecutive elements of
the ascending arrangement of all distinct elements of $\{a, b \}
\cup
 \{ \f{\ell}{n} + \vep \in (a, b) : \ell \in \bb{Z} \} \cup
 \{ \f{\ell}{n} - \vep \in (a, b) : \ell \in \bb{Z} \} $.  Then,
 $C(\se) \geq \min \{ C(\al), \; C(\ba) \}$ for any $\se \in (\al, \ba)$.
\eeL

\bpf

By Lemma \ref{constant}, both $g(\se)$ and $h(\se)$ are constants for any $\se \in (\al, \ba)$. Hence, we can drop the argument and write
$g(\se) = g, \; h(\se) = h$ and $C(\se) = S(n, g, h, \se)$.

For $\se \in (\al, \ba)$, define interval $[\al + \eta, \ba - \eta]$ with {\small $0 < \eta < \min \li ( \se - \al, \ba - \se, \f{\ba - \al} { 2
} \ri )$}. Then, \[ C(\se) \geq \min_{\mu \in [\al + \eta, \ba - \eta] } C(\mu).
\]
From the assumption that for any interval $\mscr{I}$, the probability $\Pr \li \{ Y_n \in  \mscr{I} \mid \se \ri \}$ is a continuous unimodal
function of $\se \in \Se$, we can see that, for {\small $0 < \eta < \min \li ( \se - \al, \ba - \se, \f{\ba - \al} { 2 } \ri )$}, one of the
following three cases must be true: (1) $C(\mu)$ decreases monotonically for $\mu \in [\al + \eta, \ba - \eta]$; (2) $C(\mu)$ increases
monotonically for $\mu \in [\al + \eta, \ba - \eta]$; (3) there exists a number $\se \in (\al + \eta, \ba - \eta)$ such that $C(\mu)$ increases
monotonically for $\mu \in [\al + \eta, \se]$ and decreases monotonically for $\mu \in (\se, \ba - \eta]$. It follows that
\[
C(\se) \geq \min_{\mu \in [\al + \eta, \ba - \eta]} C(\mu) = \min \{ C(\al + \eta), \; C(\ba - \eta) \}
\]
for {\small $0 < \eta < \min \li ( \se - \al, \ba - \se, \f{\ba - \al} { 2 } \ri )$}. By Lemma \ref{lem_lim}, both $\lim_{\eta \downarrow 0}
C(\al + \eta)$ and $\lim_{\eta \downarrow 0}C(\ba - \eta)$ exist and \bee C(\se) & \geq & \lim_{\eta \downarrow 0} \; \min \{
C(\al + \eta), \; C(\ba - \eta) \}\\
& = & \min \li \{ \lim_{\eta \downarrow 0} C(\al + \eta), \; \lim_{\eta \downarrow 0} C(\ba - \eta) \ri \} \geq \min \{ C(\al), \; C(\ba) \}
\eee for any $\se \in (\al, \ba)$. \epf

\bsk

Finally, to show Theorem \ref{thm_abs}, note that the statement about the coverage probability  follows
immediately from Lemma \ref{inbetween}.  The number of elements of the finite set can be calculated by using the
property of the ceiling and floor functions.

\sect{Proof of Theorem \ref{thm_rev} }

Define \[ C(\se)  =  \Pr \li \{ \li | \f{Y_n}{n} - \se \ri | < \vep \se \mid \se \ri \} = \Pr \li \{ g(\se) \leq Y_n \leq h(\se) \mid \se \ri \}
\] where
\be \la{ghrev}
 g(\se) = \lf n \se (1 - \vep) \rf + 1 , \qqu h(\se) = \lc n \se (1 + \vep ) \rc - 1.
\ee It should be noted that $C(\se), \; g(\se)$ and $h(\se)$ are actually multivariate functions of $\se, \; \vep$ and $n$.  For simplicity of
notations, we drop the arguments $n$ and $\vep$ throughout the proof of Theorem \ref{thm_rev}.

We need some preliminary results.

\beL \la{minus_rev} Let $\se_\ell = \f{\ell}{n (1 + \vep)}$ where $\ell \in \bb{Z}$. Then, $h(\se) = h(\se_{\ell + 1}) = \ell$ for any $\se \in
(\se_\ell, \se_{\ell +1})$. \eeL

\bpf For $\se \in ( \se_\ell, \; \se_{\ell + 1})$, we have $0 < n (1 + \vep) \li (\se - \se_\ell \ri ) < 1$ and \bee h (\se)  & = &
\lc n \se (1 + \vep) \rc - 1\\
& = & \lc n \se_\ell (1 + \vep) + (1 + \vep) ( \se - \se_\ell ) \rc - 1\\
& = & \li \lc n \li [ \f{\ell}{n}  + (1 + \vep) ( \se - \se_\ell)  \ri ]  \ri \rc - 1\\
& = & \ell - 1 + \li \lc n (1 +
\vep) \li (\se - \se_\ell \ri )  \ri \rc\\
& = & \ell\\
& = & \li \lc n \li [ \f{\ell + 1}{n (1 + \vep)}  \times (1 + \vep) \ri ]  \ri \rc - 1 = h(\se_{\ell + 1}). \eee

\epf

\beL \la{plus_rev} Let $\se_\ell = \f{\ell}{n (1 - \vep)}$ where $\ell \in \bb{Z}$. Then, $g(\se) = g(\se_{\ell})= \ell + 1$ for any $\se \in
(\se_\ell, \se_{\ell +1})$. \eeL

 \bpf
For $\se \in \li ( \se_\ell, \; \se_{\ell + 1} \ri )$,  we have $-1 < n (1 - \vep)  \li (\se - \se_{\ell + 1} \ri ) < 0$ and \bee
g (\se) & = &  \lf n \se (1 - \vep) \rf + 1 \\
& = & \lf n[ \se_{\ell + 1} (1 - \vep) + (1 - \vep) (\se - \se_{\ell + 1}) ] \rf + 1\\
& = &  \li \lf n \times \f{ \ell + 1 } { n (1 - \vep) } \times (1 -
\vep)  \ri \rf +
\lf n (1 - \vep) ( \se - \se_{\ell + 1} ) \rf +  1 \\
& = & \li \lf n \times \f{ \ell + 1 } { n (1 - \vep) } \times (1 -
\vep)  \ri \rf   - 1 +  1 \\
& = & \ell + 1\\
& = & \li \lf n \times \f{ \ell } { n (1 - \vep) } \times (1 - \vep) \ri \rf +  1 = g(\se_\ell). \eee

\epf

\beL \la{constant_rev} Let $\al < \ba$ be two consecutive elements
of the ascending arrangement of all distinct elements of $\{a, b
\} \cup
 \{ \f{\ell}{n (1 - \vep)}  \in (a, b) : \ell \in \bb{Z} \} \cup
 \{ \f{\ell}{n (1 + \vep)} \in (a, b) : \ell \in \bb{Z} \} $.
Then, both $g(\se)$ and $h (\se)$ are constants for any $\se \in (\al, \ba)$.
 \eeL

 \bpf
Since $\al$ and $\ba$ are two consecutive elements of the ascending arrangement of all distinct elements of the set, it must be true that there
is no integer $\ell$ such that {\small $\al < \f{\ell}{n (1 - \vep)}  < \ba$} or {\small $\al < \f{\ell}{n (1 + \vep)} < \ba$}. It follows that
there exist two integers $\ell$ and $\ell^\prime$ such that {\small $(\al, \ba) \subseteq \li ( \f{\ell}{n (1 - \vep)},  \f{\ell + 1}{n (1 -
\vep)} \ri )$} and {\small $(\al, \ba) \subseteq \li ( \f{\ell^\prime}{n (1 + \vep)}, \f{\ell^\prime + 1}{n (1 + \vep)} \ri )$}. Applying Lemma
\ref{minus_rev} and Lemma \ref{plus_rev}, we have {\small $g(\se) = g \li ( \f{\ell}{n (1 - \vep)} \ri )$} and {\small $h(\se) = h \li
(\f{\ell^\prime + 1}{n (1 + \vep)} \ri )$} for any $\se \in (\al, \ba)$.

 \epf

\beL \la{lem_lim_rev}
 For any $\se \in (0,1)$, $\lim_{\eta \downarrow 0} C(\se + \eta) \geq C(\se)$
 and $\lim_{\eta \downarrow 0} C(\se - \eta) \geq C(\se)$.
\eeL

\bpf

Observing that $h(\se + \eta) \geq h(\se)$ for any $\eta > 0$ and that
 \bee g(\se + \eta) & =  &  \lf n (\se + \eta) (1 - \vep)
\rf + 1  \\
& = &  \lf n \se (1 - \vep ) \rf + 1 + \lf n \se (1 - \vep ) -
\lf n \se ( 1  - \vep ) \rf + n \eta (1 - \vep) \rf \\
 & = & \lf n \se ( 1 - \vep )
\rf + 1   = g(\se) \eee for $0 < \eta < \f{ 1 + \lf n \se( 1 - \vep ) \rf -  n \se( 1 - \vep )} {n (1 - \vep)}$, we have \be \la{ineqa_rev} S(n,
g(\se + \eta), h (\se + \eta), \se + \eta ) \geq S(n, g(\se), h (\se), \se + \eta ) \ee for $0 < \eta < \f{ 1 + \lf n \se( 1 - \vep ) \rf - n
\se( 1 - \vep )} {n (1 - \vep)}$.  Since
 \bee h(\se + \eta) & =  &  \lc n (\se + \eta) (1 + \vep)
\rc - 1\\
& = &  \lc n \se (1 + \vep) \rc - 1 + \lc n \se (1 + \vep) - \lc n \se (1 + \vep) \rc + n \eta (1 + \vep) \rc, \eee we have {\small \[ h(\se +
\eta) = \bec \lc n \se(1 + \vep) \rc & \tx{for} \; n \se( 1 + \vep) = \lc n \se( 1  + \vep) \rc \; \tx{and} \; 0 < \eta <
\f{1}{n (1 + \vep)},\\
\lc n \se( 1 + \vep) \rc  - 1 & \tx{for} \; n \se( 1 + \vep) \neq \lc n \se( 1 + \vep) \rc \; \tx{and} \; 0 < \eta < \f{\lc n \se( 1 + \vep) \rc
- n \se( 1 + \vep)}{n (1 + \vep)}. \eec
\]}
It follows that both $g(\se + \eta)$ and $h(\se + \eta)$ are independent of $\eta$ if $\eta > 0$ is small enough. Since $S(n, g, h, \se + \eta)$
is continuous with respect to $\eta$ for fixed $g$ and $h$, we have that $\lim_{\eta \downarrow 0} S(n, g(\se + \eta), h (\se + \eta), \se +
\eta )$ exists.  As a result, \bee \lim_{\eta \downarrow 0} C(\se + \eta) & = & \lim_{\eta \downarrow 0} S(n,
g(\se + \eta), h (\se + \eta), \se + \eta )\\
& \geq & \lim_{\eta \downarrow 0} S(n, g(\se), h (\se), \se + \eta ) = S(n,  g(\se), h (\se), \se ) = C(\se), \eee where the inequality follows
from (\ref{ineqa_rev}).

Observing that $g(\se - \eta) \leq g(\se)$ for any $\eta > 0$ and that
 \bee h(\se - \eta) & =  & \lc n( \se - \eta) (1 + \vep)
\rc - 1\\
& = & \lc n \se( 1 + \vep) \rc - 1 + \lc n \se( 1 + \vep) -
\lc n \se( 1  + \vep) \rc - n \eta ( 1 + \vep) \rc\\
 & = & \lc n \se( 1 + \vep)
\rc - 1  = h(\se) \eee for $0 < \eta < \f{ 1 + n \se( 1 + \vep) - \lc n \se( 1 + \vep) \rc } {n ( 1 + \vep)}$, we have \be \la{ineqb_rev} S(n,
g(\se - \eta), h (\se - \eta), \se - \eta ) \geq S(n, g(\se), h (\se), \se - \eta ) \ee for {\small $0 < \eta < \min \li \{ \se, \f{ 1 + n \se(
1 + \vep) - \lc n \se( 1 + \vep) \rc } {n ( 1 + \vep)} \ri \}$}. Since \bee g(\se - \eta) & =  & \lf n( \se - \eta)(1 - \vep)
\rf + 1 \\
& = & \lf n \se( 1 - \vep) \rf + 1 + \lf n \se( 1 - \vep) - \lf n \se( 1 - \vep) \rf - n \eta (1 - \vep) \rf, \eee  we have {\small \[ g(\se -
\eta) = \bec \lf n \se( 1 - \vep) \rf & \tx{for} \; n \se( 1 - \vep) = \lf n \se( 1  - \vep) \rf \; \tx{and} \; 0 < \eta <
\f{1}{n (1 - \vep)},\\
\lf n \se( 1 - \vep) \rf  + 1 & \tx{for} \; n \se( 1 - \vep) \neq \lf n \se( 1 - \vep) \rf \; \tx{and} \; 0 < \eta < \f{n \se( 1 - \vep) - \lf n
\se( 1 - \vep) \rf }{n (1 - \vep)}. \eec
\]}
It follows that both $g(\se - \eta)$ and $h(\se - \eta)$ are independent of $\eta$ if $\eta > 0$ is small enough. Since $S(n, g, h, \se - \eta)$
is continuous with respect to $\eta$ for fixed $g$ and $h$, we have that $\lim_{\eta \downarrow 0} S(n, g(\se - \eta), h (\se - \eta), \se -
\eta )$ exists. Hence, \bee \lim_{\eta \downarrow 0} C(\se - \eta) & =
& \lim_{\eta \downarrow 0} S(n, g(\se - \eta), h (\se - \eta), \se -\eta )\\
& \geq & \lim_{\eta \downarrow 0} S(n, g(\se), h (\se), \se -\eta ) = S(n, g(\se), h (\se), \se ) = C(\se), \eee where the inequality follows
from (\ref{ineqb_rev}).

\epf

By a similar argument as that of Lemma \ref{inbetween}, we have \beL
\la{inbetween_rev} Let $\al < \ba$ be two consecutive elements of
the ascending arrangement of all distinct elements of $\{a, b \}
\cup
 \{ \f{\ell}{n (1 - \vep)}  \in (a, b) : \ell \in \bb{Z} \} \cup
 \{ \f{\ell}{n (1 + \vep)} \in (a, b) : \ell \in \bb{Z} \} $.  Then,
 $C(\se) \geq \min \{ C(\al), \; C(\ba) \}$ for any $\se \in (\al, \ba)$.
\eeL

\bsk

Finally, to show Theorem \ref{thm_rev}, note that the statement about the coverage probability  follows
immediately from Lemma \ref{inbetween_rev}.  The number of elements of the finite set can be calculated by using
the property of the ceiling and floor functions.

\sect{Proof of Theorem \ref{Chensabs}}  \la{Chensabs_app}

We need some preliminary results.

\beL \la{lemvipaabs} Assume that $a + \vep \leq b - \vep$.  Then, \bel &  &
 \{ | \wt{\bs{\se}}_n - \se | \geq \vep   \} = \{ |
\wh{\bs{\se}}_n - \se | \geq \vep  \} \qqu \tx{for $\se \in \li [ a + \vep, b - \vep \ri ]$}, \la{ref88aabs}\\
&  & \{ | \wt{\bs{\se}}_n - \se | \geq \vep  \} = \{ \wh{\bs{\se}}_n  \geq \se + \vep \} \qqu \tx{for $\se \in \li [ a, a +
\vep \ri )$}, \la{ref88babs} \\
&  & \{ | \wt{\bs{\se}}_n - \se | \geq \vep \} = \{ \wh{\bs{\se}}_n  \leq \se - \vep \} \qqu \tx{for $\se \in \li ( b - \vep, b \ri ]$}.
\la{ref88cabs}
 \eel

\eeL

\bpf  To prove (\ref{ref88aabs}), recalling the definition of the range-preserving estimator, we have that \be \la{comb33881abs}
 \{ | \wt{\bs{\se}}_n - \se | \geq \vep,  \; \wh{\bs{\se}}_n \in [a, b]  \} = \{ |
\wh{\bs{\se}}_n - \se | \geq \vep, \; \wh{\bs{\se}}_n \in [a, b] \} \ee for $\se \in \Se$.  For $\se \in \li [ a + \vep, b - \vep \ri ]$, we
have $\se + \vep > a$ and
\[
 \{ \wt{\bs{\se}}_n \geq \se + \vep,  \; \wh{\bs{\se}}_n < a \}   \subseteq  \{ \wt{\bs{\se}}_n > a,  \; \wh{\bs{\se}}_n < a \}  = \emptyset, \]
 \[
 \{ \wh{\bs{\se}}_n \geq \se + \vep,  \; \wh{\bs{\se}}_n < a \}   \subseteq  \{ \wh{\bs{\se}}_n > a,  \; \wh{\bs{\se}}_n < a \}  = \emptyset, \]
 which implies that
 \be
 \la{comb33882abs}
 \{ \wt{\bs{\se}}_n \geq \se + \vep,  \; \wh{\bs{\se}}_n < a \}   =  \{ \wh{\bs{\se}}_n \geq  \se + \vep,  \; \wh{\bs{\se}}_n < a \}
 \ee
for $\se \in \li [ a + \vep, b - \vep \ri ]$.  On the other hand, for $\se \in \li [ a + \vep, b - \vep \ri ]$,  we have $ \se - \vep < b$ and
\[
 \{ \wt{\bs{\se}}_n \leq \se - \vep,  \; \wh{\bs{\se}}_n > b \}   \subseteq  \{ \wt{\bs{\se}}_n < b,  \; \wh{\bs{\se}}_n > b \}  = \emptyset, \]
 \[
 \{ \wh{\bs{\se}}_n \leq \se - \vep,  \; \wh{\bs{\se}}_n > b \}   \subseteq  \{ \wh{\bs{\se}}_n < b,  \; \wh{\bs{\se}}_n > b \}  = \emptyset, \]
 which implies that
 \be
 \la{comb33883abs}
 \{ \wt{\bs{\se}}_n \leq \se - \vep,  \; \wh{\bs{\se}}_n > b \}   =  \{ \wh{\bs{\se}}_n \leq \se - \vep,  \; \wh{\bs{\se}}_n > b \}
 \ee
for $\se \in \li [ a + \vep, b - \vep \ri ]$.  Combining (\ref{comb33881abs}), (\ref{comb33882abs}), and (\ref{comb33883abs}) completes the
proof of (\ref{ref88aabs}).

To prove (\ref{ref88babs}), note that, for $\se \in \li [ a, a + \vep \ri )$,  we have $\se - \vep < a$ and
\[
 \{ \wt{\bs{\se}}_n \leq \se - \vep \}   \subseteq  \{ \wt{\bs{\se}}_n < a \}  = \emptyset. \]
 On the other hand, for $\se \in \li [ a, a + \vep \ri )$,  we have $a <  \se + \vep < b$ and
\[
 \{ \wt{\bs{\se}}_n \geq \se + \vep \}   = \{ \wh{\bs{\se}}_n \geq \se + \vep \}. \]  This proves (\ref{ref88babs}).

To prove (\ref{ref88cabs}), note that, for $\se \in \li ( b - \vep, b \ri ]$,  we have $\se + \vep> b$ and
\[
 \{ \wt{\bs{\se}}_n \geq \se + \vep \}   \subseteq  \{ \wt{\bs{\se}}_n > b \}  = \emptyset. \]
 On the other hand, for $\se \in \li ( b - \vep, b \ri ]$,  we have $a < \se - \vep < b$ and
\[
 \{ \wt{\bs{\se}}_n \leq \se - \vep \}   = \{ \wh{\bs{\se}}_n \leq  \se - \vep \}. \]  This proves (\ref{ref88cabs}).   The proof of the
 lemma is thus completed.

\epf

\beL  \la{comb33886abs}

Assume that $a + \vep > b - \vep$.  Then, \bel &  &  \{ | \wt{\bs{\se}}_n - \se | \geq \vep  \} = \emptyset  \qqu \tx{for
$\se \in \li ( b - \vep,  a + \vep \ri )$}, \la{89933aabs}\\
&  & \{ | \wt{\bs{\se}}_n - \se | \geq \vep \} = \{ \wh{\bs{\se}}_n  \geq \se + \vep \} \qqu \tx{for $\se \in \li [ a, b -
\vep \ri ]$}, \la{89933babs}\\
&  & \{ | \wt{\bs{\se}}_n - \se | \geq \vep \} = \{ \wh{\bs{\se}}_n  \leq \se - \vep \} \qqu \tx{for $\se \in \li [ a + \vep, b \ri ]$}.
\la{89933cabs} \eel

\eeL

\bpf

To show (\ref{89933aabs}), note that, for $\se \in \li ( b - \vep, a + \vep \ri )$,  we have $\se - \vep < a$ and
\[
 \{ \wt{\bs{\se}}_n \leq \se - \vep \}   \subseteq  \{ \wt{\bs{\se}}_n < a \}  = \emptyset. \]
 On the other hand, for $\se \in \li ( b - \vep, a + \vep \ri )$,  we have $\se  + \vep > b$ and
\[
 \{ \wt{\bs{\se}}_n \geq \se  + \vep \}   \subseteq  \{ \wt{\bs{\se}}_n > b \}  = \emptyset. \]  This proves (\ref{89933aabs}).

To show (\ref{89933babs}), note that for $\se \in \li [ a, b - \vep \ri ]$,  we have $\se  - \vep < a$ and
\[
 \{ \wt{\bs{\se}}_n \leq \se - \vep \}   \subseteq  \{ \wt{\bs{\se}}_n < a \}  = \emptyset. \]
 On the other hand, for $\se \in \li [ a, b - \vep \ri ]$,  we have $a < \se + \vep \leq b$ and
\[
 \{ \wt{\bs{\se}}_n \geq \se + \vep \}   = \{ \wh{\bs{\se}}_n \geq \se + \vep \}. \]  This proves (\ref{89933babs}).

To show (\ref{89933cabs}), note that for $\se \in \li [ a + \vep, b \ri ]$,  we have $\se + \vep > b$ and
\[
 \{ \wt{\bs{\se}}_n \geq \se + \vep \}   \subseteq  \{ \wt{\bs{\se}}_n > b \}  = \emptyset. \]
 On the other hand, for $\se \in \li [ a + \vep, b \ri ]$,  we have $a \leq \se - \vep < b$ and
\[
 \{ \wt{\bs{\se}}_n \leq \se - \vep \}   = \{ \wh{\bs{\se}}_n \leq \se - \vep \}. \]  This proves (\ref{89933cabs}).  The proof of the
 lemma is thus completed.

\epf

\beL

\la{lem3388a96abs}

Assume that $a + \vep \leq b - \vep$.  Then, the minimum of $\Pr \{ | \wt{\bs{\se}}_n - \se | < \vep  \mid \se \}$ with respect to $\se \in \li
[ a + \vep, b - \vep \ri ]$ is attained at the finite set {\small \[ A_0 =  \li \{ a + \vep, b - \vep \ri \} \bigcup \li \{ \f{k}{n} - \vep \in
\li ( a + \vep, b - \vep \ri ): k \in \bb{Z} \ri \} \bigcup \li \{  \f{k}{n} + \vep \in \li ( a + \vep, b - \vep \ri ): k \in \bb{Z} \ri \}
\]
} \eeL

\bpf

From $\se \in \li [ a + \vep, b - \vep \ri ]$, it follows from (\ref{ref88aabs}) of Lemma \ref{lemvipaabs} that
\[
\Pr \{ | \wt{\bs{\se}}_n - \se | < \vep \mid \se \} = \Pr \{ | \wh{\bs{\se}}_n - \se | < \vep \mid \se \}.
\]
Hence, by Theorem \ref{thm_abs}, the lemma follows.

\epf

\beL \la{lem3388b96abs} Assume that $a + \vep \leq b - \vep$.   Then, the minimum of $\Pr \{ | \wt{\bs{\se}}_n - \se | < \vep \mid \se \}$ with
respect to $\se \in \li [a, a + \vep \ri ]$ is attained at the finite set {\small \[ A_1 =  \li \{ a,  a + \vep \ri \} \bigcup \li \{ \f{k}{n} -
\vep \in \li (a, a + \vep \ri ): k \in \bb{Z} \ri \}.
\]
} \eeL

\bpf

For $\se \in \li [a, a + \vep \ri )$, it follows from (\ref{ref88babs}) of Lemma \ref{lemvipaabs} that $\Pr \{ | \wt{\bs{\se}}_n - \se | < \vep
\mid \se \} = \Pr \{ \wh{\bs{\se}}_n < \se + \vep \mid \se \}$.  Let $\al < \ba$ be two consecutive elements of $A_1$.  It follows from Lemma
\ref{minus} that
\[
\Pr \{ | \wt{\bs{\se}}_n - \se | < \vep \mid \se \} = \Pr \{ \wh{\bs{\se}}_n < \se + \vep \mid \se \} = \Pr \{ \wh{\bs{\se}}_n < \se^\star +
\vep \mid \se \}
\]
for $\se \in (\al, \ba)$, where $\se^\star = \f{\al + \ba}{2}$.  Now let $\eta \in (0, \f{\ba - \al}{2} )$.  Recalling the assumption that for
any interval $\mscr{I}$, the probability $\Pr \{ Y_n \in \mscr{I} \mid \se \}$ is a unimodal function of $\se \in [a, b]$, we have that {\small
\bee \Pr \{ \wh{\bs{\se}}_n < \se + \vep \mid \se \} & \geq & \min \li [ \Pr \{ \wh{\bs{\se}}_n < (\al + \eta) + \vep \mid \al + \eta \}, \; \Pr
\{ \wh{\bs{\se}}_n < (\ba - \eta)  + \vep \mid \ba - \eta \} \ri ]\\
& =  & \min \li [  \Pr \{ Y_n \leq h(\al + \eta) \mid \al + \eta \}, \; \Pr \{ Y_n \leq h(\ba - \eta) \mid \ba - \eta \} \ri ] \eee} for any
$\se \in (\al, \ba)$, where the function $h(.)$ is defined by (\ref{ghabs}).  From Lemma \ref{minus}, we know that
\[
h(\al + \eta) = h(\ba - \eta) = h(\ba).
\]
Hence, \bee \Pr \{ \wh{\bs{\se}}_n < \se + \vep \mid \se  \}   \geq    \min \li [  \Pr \{ Y_n \leq h(\ba) \mid \al + \eta \}, \; \Pr \{ Y_n \leq
h(\ba) \mid \ba - \eta \} \ri ] \eee for any $\se \in (\al, \ba)$.  Recalling the assumption that for any interval $\mscr{I}$, the probability
$\Pr \{ Y_n \in \mscr{I} \mid \se \}$ is a continuous function of $\se \in [a, b]$, we have \bee \Pr \{ | \wt{\bs{\se}}_n - \se | < \vep \mid
\se \} & = & \Pr \{ \wh{\bs{\se}}_n < \se + \vep \mid \se \}\\
& \geq &  \lim_{\eta \downarrow 0} \min \li [  \Pr \{ Y_n \leq h(\ba) \mid \al + \eta \}, \; \Pr \{ Y_n \leq
h(\ba) \mid \ba - \eta \} \ri ]\\
&  = &   \min \li [  \lim_{\eta \downarrow 0} \Pr \{ Y_n \leq h(\ba) \mid \al + \eta \}, \; \lim_{\eta \downarrow 0} \Pr \{ Y_n \leq h(\ba) \mid
\ba - \eta \} \ri ]\\
&  = &   \min \li [  \Pr \{ Y_n \leq h(\ba) \mid \al  \}, \; \Pr \{ Y_n \leq h(\ba) \mid \ba  \} \ri ]\\
& \geq &  \min \li [  \Pr \{ Y_n \leq h(\al) \mid \al  \}, \; \Pr \{ Y_n \leq h(\ba) \mid \ba \} \ri ]\\
& = & \min \li [  \Pr \{ \wh{\bs{\se}}_n < \al + \vep \mid \al \}, \; \Pr \{ \wh{\bs{\se}}_n < \ba + \vep \mid \ba \}  \ri ]\\
& = & \min \li [  \Pr \{ | \wt{\bs{\se}}_n - \al | < \vep \mid \al \}, \; \Pr \{ | \wt{\bs{\se}}_n - \ba | < \vep \mid \ba \}  \ri ] \eee  for
any $\se \in (\al, \ba)$. This immediately leads to the result of the lemma. \epf

\beL \la{lem3388c96abs} Assume that $a + \vep \leq b - \vep$.  Then, the minimum of $\Pr \{ | \wt{\bs{\se}}_n - \se | < \vep \mid \se \}$ with
respect to $\se \in \li [ b - \vep, b \ri ]$ is attained at the finite set {\small \[ A_2 =  \li \{ b, b - \vep \ri \} \bigcup \li \{ \f{k}{n} +
\vep \in \li ( b - \vep, b \ri ): k \in \bb{Z} \ri \}.
\]
} \eeL

\bpf

For $\se \in \li ( b - \vep, b \ri ]$, it follows from (\ref{ref88cabs}) of Lemma \ref{lemvipaabs} that $\Pr \{ | \wt{\bs{\se}}_n - \se | < \vep
\mid \se \} = \Pr \{ \wh{\bs{\se}}_n > \se - \vep \mid \se \}$.  Let $\al < \ba$ be two consecutive elements of $A_2$.  It follows from Lemma
\ref{plus} that
\[
\Pr \{ | \wt{\bs{\se}}_n - \se | < \vep \mid \se \} = \Pr \{ \wh{\bs{\se}}_n > \se - \vep \mid \se \} = \Pr \{ \wh{\bs{\se}}_n > \se^* - \vep
\mid \se \}
\]
for $\se \in (\al, \ba)$, where $\se^* = \f{\al + \ba}{2}$.  Now let $\eta \in (0, \f{\ba - \al}{2} )$.  Invoking the assumption that for any
interval $\mscr{I}$, the probability $\Pr \{ Y_n \in \mscr{I} \mid \se \}$ is a unimodal function of $\se \in [a, b]$, we have {\small \bee \Pr
\{ \wh{\bs{\se}}_n
> \se - \vep \mid \se  \} &  \geq &  \min \li [  \Pr \{ \wh{\bs{\se}}_n > (\al + \eta) - \vep \mid \al + \eta \}, \; \Pr \{
\wh{\bs{\se}}_n > (\ba - \eta) - \vep  \mid \ba - \eta \} \ri ]\\
& =  & \min \li [  \Pr \{ Y_n \geq g(\al + \eta) \mid \al + \eta \}, \; \Pr \{ Y_n \geq g(\ba - \eta) \mid \ba - \eta \} \ri ] \eee} for any
$\se \in (\al, \ba)$, where the function $g(.)$ is defined by (\ref{ghabs}).  From Lemma \ref{plus}, we know that
\[
g(\al + \eta) = g(\ba - \eta) = g(\al).
\]
Hence, \bee \Pr \{ \wh{\bs{\se}}_n > \se - \vep \mid \se  \}  \geq   \min \li [  \Pr \{ Y_n \geq g(\al) \mid \al + \eta \}, \; \Pr \{ Y_n \geq
g(\al) \mid \ba - \eta \} \ri ] \eee for any $\se \in (\al, \ba)$.  Recalling the assumption that for any interval $\mscr{I}$, the probability
$\Pr \{ Y_n \in \mscr{I} \mid \se \}$ is a continuous function of $\se \in [a, b]$, we have \bee \Pr \{ | \wt{\bs{\se}}_n - \se | < \vep \mid
\se \} & = & \Pr \{ \wh{\bs{\se}}_n > \se - \vep \mid \se \}\\
 & \geq &  \lim_{\eta \downarrow 0} \min \li [  \Pr \{ Y_n \geq g(\al) \mid \al +
\eta \}, \; \Pr \{ Y_n \geq
g(\al) \mid \ba - \eta \} \ri ]\\
&  = &   \min \li [  \lim_{\eta \downarrow 0} \Pr \{ Y_n \geq g (\al) \mid \al + \eta \}, \; \lim_{\eta \downarrow 0} \Pr \{ Y_n \geq g(\al)
\mid \ba - \eta \} \ri ]\\
&  = &   \min \li [  \Pr \{ Y_n \geq g(\al) \mid \al  \}, \; \Pr \{ Y_n \geq g(\al) \mid \ba  \} \ri ]\\
& \geq &  \min \li [  \Pr \{ Y_n \geq g(\al) \mid \al  \}, \; \Pr \{ Y_n \geq g(\ba) \mid \ba \} \ri ]\\
& = &  \min \li [  \Pr \{ \wh{\bs{\se}}_n > \al - \vep \mid \al \}, \; \Pr \{ \wh{\bs{\se}}_n > \ba - \vep \mid \ba \}  \ri ]\\
& = &  \min \li [  \Pr \{ | \wt{\bs{\se}}_n - \al | < \vep \mid \al \}, \;  \Pr \{ | \wt{\bs{\se}}_n - \ba | < \vep \mid \ba \} \ri ]
  \eee for any $\se \in (\al, \ba)$. This immediately leads to the result of the lemma.  \epf

By virtue of (\ref{89933aabs}) of Lemma \ref{comb33886abs}, we have the following result.

\beL

\la{lem3388a9abs}

Assume that $a + \vep > b - \vep$.  Then, $\Pr \{ | \wt{\bs{\se}}_n - \se | < \vep \mid \se \} = 1$ for $\se \in \li ( b - \vep, a + \vep \ri
)$.  \eeL

By virtue of (\ref{89933babs}) of Lemma \ref{comb33886abs} and a similar argument as that of Lemma \ref{lem3388b96abs}, we have the following
result.

 \beL \la{lem3388b9abs} Assume that $b \geq  a + \vep > b - \vep$.   Then, the minimum of $\Pr \{ | \wt{\bs{\se}}_n - \se | < \vep \mid \se \}$ with respect to $\se \in \li [a, b - \vep \ri ]$
is attained at the finite set {\small \[ A_3 =  \li \{ a,  b - \vep \ri \} \bigcup \li \{ \f{k}{n} - \vep \in \li (a,  b - \vep \ri ): k \in
\bb{Z} \ri \}.
\]
} \eeL

By virtue of (\ref{89933cabs}) of Lemma \ref{comb33886abs} and a similar argument as that of Lemma \ref{lem3388c96abs}, we have the following
result.

\beL \la{lem3388c9abs} Assume that $b \geq  a + \vep > b - \vep$.   Then, the minimum of $\Pr \{ | \wt{\bs{\se}}_n - \se | < \vep \mid \se \}$
with respect to $\se \in \li [ a + \vep, b \ri ]$ is attained at the finite set {\small \[ A_4 =  \li \{ b,  a + \vep \ri \} \bigcup \li \{
\f{k}{n} + \vep \in \li ( a + \vep, b \ri ): k \in \bb{Z} \ri \}.
\]
} \eeL

We are now in a position to prove the theorem.  We need to consider three cases as follows.

In the case of $a + \vep \leq b - \vep$, it follows from Lemmas \ref{lem3388a96abs}--\ref{lem3388c96abs} that, the minimum of $\Pr \{ |
\wt{\bs{\se}}_n - \se | < \vep \mid \se \}$ with respect to $\se \in [a, b]$ is attained at the finite set $( A_0 \cup A_1 \cup A_2 ) \cap [a,
b] = A \cap [a, b]$, where $A$ is the set defined by (\ref{Aabs}).

In the case of $b \geq a + \vep
> b - \vep$, it follows from Lemmas \ref{lem3388a9abs}--\ref{lem3388c9abs} that, the minimum of $\Pr \{ | \wt{\bs{\se}}_n - \se | < \vep \mid
\se \}$ with respect to $\se \in [a, b]$ is attained at the finite set $(A_3 \cup A_4) \cap [a, b] = A \cap [a, b]$.

In the case of $a + \vep > b$, we have $b - \vep < a < b < a + \vep$.  It follows from Lemma \ref{lem3388a9abs} that $\Pr \{ | \wt{\bs{\se}}_n -
\se | < \vep \mid \se \} = 1$ for $\se \in [a, b]$.  Thus, the minimum of $\Pr \{ | \wt{\bs{\se}}_n - \se | < \vep \mid \se \}$ with respect to
$\se \in [a, b]$ is equal to $1$, which is  attained at $\{a, b \} = A \cap [a, b]$.

The number of elements of the finite set can be calculated by using the property of the ceiling and floor functions. This completes the proof of
the theorem.

\sect{Proof of Theorem \ref{Chens}} \la{Chens_app}

We need some preliminary results.  The following Lemmas \ref{lemvipa}--\ref{lem3388a96}  are more general but similar to the results of
\cite{Gamrot}. To justify these results, we also follow similar arguments as that of \cite{Gamrot}.

\beL \la{lemvipa}
Assume that $\f{a}{1 - \vep} \leq \f{b}{1 + \vep} $.  Then, \bel &  &
 \{ | \wt{\bs{\se}}_n - \se | \geq \vep \; \se  \} = \{ |
\wh{\bs{\se}}_n - \se | \geq \vep \; \se \} \qqu \tx{for $\se \in \li [ \f{a}{1 - \vep}, \f{b}{1 + \vep} \ri ]$}, \la{ref88a}\\
&  & \{ | \wt{\bs{\se}}_n - \se | \geq \vep \; \se  \} = \{ \wh{\bs{\se}}_n  \geq (1 + \vep) \; \se \} \qqu \tx{for $\se \in \li [ a, \f{a}{1 -
\vep} \ri )$}, \la{ref88b} \\
&  & \{ | \wt{\bs{\se}}_n - \se | \geq \vep \; \se  \} = \{ \wh{\bs{\se}}_n  \leq (1 - \vep) \; \se \} \qqu \tx{for $\se \in \li ( \f{b}{1 +
\vep}, b \ri ]$}.  \la{ref88c}
 \eel

\eeL

\bpf  To prove (\ref{ref88a}), recalling the definition of the range-preserving estimator, we have that \be \la{comb33881}
 \{ | \wt{\bs{\se}}_n - \se | \geq \vep \; \se,  \; \wh{\bs{\se}}_n \in [a, b]  \} = \{ |
\wh{\bs{\se}}_n - \se | \geq \vep \; \se, \; \wh{\bs{\se}}_n \in [a, b] \} \ee for $\se \in \Se$.  For $\se \in \li [ \f{a}{1 - \vep}, \f{b}{1 +
\vep} \ri ]$, we have $(1 + \vep) \se > a$ and
\[
 \{ \wt{\bs{\se}}_n \geq (1 + \vep)  \se,  \; \wh{\bs{\se}}_n < a \}   \subseteq  \{ \wt{\bs{\se}}_n > a,  \; \wh{\bs{\se}}_n < a \}  = \emptyset, \]
 \[
 \{ \wh{\bs{\se}}_n \geq (1 + \vep)  \se,  \; \wh{\bs{\se}}_n < a \}   \subseteq  \{ \wh{\bs{\se}}_n > a,  \; \wh{\bs{\se}}_n < a \}  = \emptyset, \]
 which implies that
 \be
 \la{comb33882}
 \{ \wt{\bs{\se}}_n \geq (1 + \vep)  \se,  \; \wh{\bs{\se}}_n < a \}   =  \{ \wh{\bs{\se}}_n \geq (1 + \vep)  \se,  \; \wh{\bs{\se}}_n < a \}
 \ee
for $\se \in \li [ \f{a}{1 - \vep}, \f{b}{1 + \vep} \ri ]$.  On the other hand, for $\se \in \li [ \f{a}{1 - \vep}, \f{b}{1 + \vep} \ri ]$,  we
have $(1 - \vep) \se < b$ and
\[
 \{ \wt{\bs{\se}}_n \leq (1 - \vep)  \se,  \; \wh{\bs{\se}}_n > b \}   \subseteq  \{ \wt{\bs{\se}}_n < b,  \; \wh{\bs{\se}}_n > b \}  = \emptyset, \]
 \[
 \{ \wh{\bs{\se}}_n \leq (1 - \vep)  \se,  \; \wh{\bs{\se}}_n > b \}   \subseteq  \{ \wh{\bs{\se}}_n < b,  \; \wh{\bs{\se}}_n > b \}  = \emptyset, \]
 which implies that
 \be
 \la{comb33883}
 \{ \wt{\bs{\se}}_n \leq (1 - \vep)  \se,  \; \wh{\bs{\se}}_n > b \}   =  \{ \wh{\bs{\se}}_n \leq (1 - \vep)  \se,  \; \wh{\bs{\se}}_n > b \}
 \ee
for $\se \in \li [ \f{a}{1 - \vep}, \f{b}{1 + \vep} \ri ]$.  Combining (\ref{comb33881}), (\ref{comb33882}), and (\ref{comb33883}) completes the
proof of (\ref{ref88a}).

To prove (\ref{ref88b}), note that, for $\se \in \li [ a, \f{a}{1 - \vep} \ri )$,  we have $(1 - \vep) \se < a$ and
\[
 \{ \wt{\bs{\se}}_n \leq (1 - \vep)  \se \}   \subseteq  \{ \wt{\bs{\se}}_n < a \}  = \emptyset. \]
 On the other hand, for $\se \in \li [ a, \f{a}{1 - \vep} \ri )$,  we have $a < (1 + \vep) \se < b$ and
\[
 \{ \wt{\bs{\se}}_n \geq (1 + \vep)  \se \}   = \{ \wh{\bs{\se}}_n \geq (1 + \vep)  \se \}. \]  This proves (\ref{ref88b}).

To prove (\ref{ref88c}), note that, for $\se \in \li ( \f{b}{1 + \vep}, b \ri ]$,  we have $(1 + \vep) \se > b$ and
\[
 \{ \wt{\bs{\se}}_n \geq (1 + \vep)  \se \}   \subseteq  \{ \wt{\bs{\se}}_n > b \}  = \emptyset. \]
 On the other hand, for $\se \in \li ( \f{b}{1 + \vep}, b \ri ]$,  we have $a < (1 - \vep) \se < b$ and
\[
 \{ \wt{\bs{\se}}_n \leq (1 - \vep)  \se \}   = \{ \wh{\bs{\se}}_n \leq (1 - \vep)  \se \}. \]  This proves (\ref{ref88c}).   The proof of the
 lemma is thus completed.

\epf

\beL  \la{comb33886}

Assume that $\f{a}{1 - \vep} > \f{b}{1 + \vep} $.  Then, \bel &  &  \{ | \wt{\bs{\se}}_n - \se | \geq \vep \; \se  \} = \emptyset  \qqu \tx{for
$\se \in \li ( \f{b}{1 + \vep},  \f{a}{1 - \vep} \ri )$}, \la{89933a}\\
&  & \{ | \wt{\bs{\se}}_n - \se | \geq \vep \; \se  \} = \{ \wh{\bs{\se}}_n  \geq (1 + \vep) \; \se \} \qqu \tx{for $\se \in \li [ a, \f{b}{1 +
\vep} \ri ]$}, \la{89933b}\\
&  & \{ | \wt{\bs{\se}}_n - \se | \geq \vep \; \se  \} = \{ \wh{\bs{\se}}_n  \leq (1 - \vep) \; \se \} \qqu \tx{for $\se \in \li [ \f{a}{1 -
\vep}, b \ri ]$}.  \la{89933c} \eel

\eeL

\bpf

To show (\ref{89933a}), note that, for $\se \in \li ( \f{b}{1 + \vep}, \f{a}{1 - \vep} \ri )$,  we have $(1 - \vep) \se < a$ and
\[
 \{ \wt{\bs{\se}}_n \leq (1 - \vep)  \se \}   \subseteq  \{ \wt{\bs{\se}}_n < a \}  = \emptyset. \]
 On the other hand, for $\se \in \li ( \f{b}{1 + \vep}, \f{a}{1 - \vep} \ri )$,  we have $(1 + \vep) \se > b$ and
\[
 \{ \wt{\bs{\se}}_n \geq (1 + \vep)  \se \}   \subseteq  \{ \wt{\bs{\se}}_n > b \}  = \emptyset. \]  This proves (\ref{89933a}).

To show (\ref{89933b}), note that for $\se \in \li [ a, \f{b}{1 + \vep} \ri ]$,  we have $(1 - \vep) \se < a$ and
\[
 \{ \wt{\bs{\se}}_n \leq (1 - \vep)  \se \}   \subseteq  \{ \wt{\bs{\se}}_n < a \}  = \emptyset. \]
 On the other hand, for $\se \in \li [ a, \f{b}{1 + \vep} \ri ]$,  we have $a < (1 + \vep) \se \leq b$ and
\[
 \{ \wt{\bs{\se}}_n \geq (1 + \vep)  \se \}   = \{ \wh{\bs{\se}}_n \geq (1 + \vep)  \se \}. \]  This proves (\ref{89933b}).

To show (\ref{89933c}), note that for $\se \in \li [ \f{a}{1 - \vep}, b \ri ]$,  we have $(1 + \vep) \se > b$ and
\[
 \{ \wt{\bs{\se}}_n \geq (1 + \vep)  \se \}   \subseteq  \{ \wt{\bs{\se}}_n > b \}  = \emptyset. \]
 On the other hand, for $\se \in \li [ \f{a}{1 - \vep}, b \ri ]$,  we have $a \leq (1 - \vep) \se < b$ and
\[
 \{ \wt{\bs{\se}}_n \leq (1 - \vep)  \se \}   = \{ \wh{\bs{\se}}_n \leq (1 - \vep)  \se \}. \]  This proves (\ref{89933c}).  The proof of the
 lemma is thus completed.

\epf

\beL

\la{lem3388a96}

Assume that $\f{a}{1 - \vep} \leq \f{b}{1 + \vep}$.  Then, the minimum of $\Pr \{ | \wt{\bs{\se}}_n - \se | < \vep \se \mid \se \}$ with respect
to $\se \in \li [ \f{a}{1 - \vep}, \f{b}{1 + \vep} \ri ]$ is attained at the finite set {\small \[ A_0 =  \li \{ \f{a}{1 - \vep}, \f{b}{1 +
\vep} \ri \} \bigcup \li \{  \f{k}{n (1 + \vep)} \in \li ( \f{a}{1 - \vep}, \f{b}{1 + \vep} \ri ): k \in \bb{Z} \ri \} \bigcup \li \{  \f{k}{n
(1 - \vep)} \in \li ( \f{a}{1 - \vep}, \f{b}{1 + \vep} \ri ): k \in \bb{Z} \ri \}
\]
} \eeL

\bpf

From $\se \in \li [ \f{a}{1 - \vep}, \f{b}{1 + \vep} \ri ]$, it follows from (\ref{ref88a}) of Lemma \ref{lemvipa} that
\[
\Pr \{ | \wt{\bs{\se}}_n - \se | < \vep \se \mid \se \} = \Pr \{ | \wh{\bs{\se}}_n - \se | < \vep \se \mid \se \}.
\]
Hence, by Theorem \ref{thm_rev}, the lemma follows.

\epf

The following Lemmas \ref{lem3388b96} and \ref{lem3388c96}  are more general but similar to the results of \cite{Gamrot}. To justify these
results, an analysis of discontinuity is necessary.

\beL \la{lem3388b96} Assume that $\f{a}{1 - \vep} \leq \f{b}{1 + \vep}$.   Then, the minimum of $\Pr \{ | \wt{\bs{\se}}_n - \se | < \vep \se
\mid \se \}$ with respect to $\se \in \li [a, \f{a}{1 - \vep} \ri ]$ is attained at the finite set {\small \[ A_1 =  \li \{ a,  \f{a}{1 - \vep}
\ri \} \bigcup \li \{ \f{k}{n (1 + \vep)} \in \li (a, \f{a}{1 - \vep} \ri ): k \in \bb{Z} \ri \}.
\]
} \eeL

\bpf

For $\se \in \li [a, \f{a}{1 - \vep} \ri )$, it follows from (\ref{ref88b}) of Lemma \ref{lemvipa} that $\Pr \{ | \wt{\bs{\se}}_n - \se | < \vep
\se \mid \se \} = \Pr \{ \wh{\bs{\se}}_n < (1 + \vep) \se \mid \se \}$.  Let $\al < \ba$ be two consecutive elements of $A_1$.  It follows from
Lemma \ref{minus_rev} that
\[
\Pr \{ | \wt{\bs{\se}}_n - \se | < \vep \se \mid \se \} = \Pr \{ \wh{\bs{\se}}_n < (1 + \vep) \se \mid \se \} = \Pr \{ \wh{\bs{\se}}_n < (1 +
\vep) \se^\star \mid \se \}
\]
for $\se \in (\al, \ba)$, where $\se^\star = \f{\al + \ba}{2}$.  Now let $\eta \in (0, \f{\ba - \al}{2} )$.  Recalling the assumption that for
any interval $\mscr{I}$, the probability $\Pr \{ Y_n \in \mscr{I} \mid \se \}$ is a unimodal function of $\se \in [a, b]$, we have that {\small
\bee \Pr \{ \wh{\bs{\se}}_n < (1 + \vep) \se \mid \se \} & \geq & \min \li [ \Pr \{ \wh{\bs{\se}}_n < (1 + \vep) (\al + \eta) \mid \al + \eta
\}, \; \Pr \{
\wh{\bs{\se}}_n < (1 + \vep) (\ba - \eta) \mid \ba - \eta \} \ri ]\\
& =  & \min \li [  \Pr \{ Y_n \leq h(\al + \eta) \mid \al + \eta \}, \; \Pr \{ Y_n \leq h(\ba - \eta) \mid \ba - \eta \} \ri ] \eee} for any
$\se \in (\al, \ba)$, where the function $h(.)$ is defined by (\ref{ghrev}).   From Lemma \ref{minus_rev}, we know that
\[
h(\al + \eta) = h(\ba - \eta) = h(\ba).
\]
Hence, \bee \ \Pr \{ \wh{\bs{\se}}_n < (1 + \vep) \se \mid \se  \}  \geq  \min \li [ \Pr \{ Y_n \leq h(\ba) \mid \al + \eta \}, \; \Pr \{ Y_n
\leq h(\ba) \mid \ba - \eta \} \ri ] \eee for any $\se \in (\al, \ba)$. Recalling the assumption that for any interval $\mscr{I}$, the
probability $\Pr \{ Y_n \in \mscr{I} \mid \se \}$ is a continuous function of $\se \in [a, b]$, we have \bee \Pr \{ | \wt{\bs{\se}}_n - \se | <
\vep \se \mid \se \} & = & \Pr \{ \wh{\bs{\se}}_n < (1 + \vep) \se \mid \se \}\\
 & \geq &  \lim_{\eta \downarrow 0} \min \li [ \Pr \{ Y_n \leq
h(\ba) \mid \al + \eta \}, \; \Pr \{ Y_n \leq
h(\ba) \mid \ba - \eta \} \ri ]\\
&  = &   \min \li [  \lim_{\eta \downarrow 0} \Pr \{ Y_n \leq h(\ba) \mid \al + \eta \}, \; \lim_{\eta \downarrow 0} \Pr \{ Y_n \leq h(\ba) \mid
\ba - \eta \} \ri ]\\
&  = &   \min \li [  \Pr \{ Y_n \leq h(\ba) \mid \al  \}, \; \Pr \{ Y_n \leq h(\ba) \mid \ba  \} \ri ]\\
& \geq &  \min \li [  \Pr \{ Y_n \leq h(\al) \mid \al  \}, \; \Pr \{ Y_n \leq h(\ba) \mid \ba \} \ri ]\\
& = &   \min \li [  \Pr \{ \wh{\bs{\se}}_n < (1 + \vep) \al \mid \al \}, \; \Pr \{ \wh{\bs{\se}}_n < (1 + \vep) \ba \mid \ba \} \ri ]\\
& = &  \min \li [ \Pr \{ | \wt{\bs{\se}}_n - \al | < \vep \al \mid \al \}, \; \Pr \{ | \wt{\bs{\se}}_n - \ba | < \vep \ba \mid \ba \}  \ri ]
\eee  for any $\se \in (\al, \ba)$. This immediately leads to the result of the lemma. \epf

\beL \la{lem3388c96} Assume that $\f{a}{1 - \vep} \leq \f{b}{1 + \vep}$.  Then, the minimum of $\Pr \{ | \wt{\bs{\se}}_n - \se | < \vep \se \mid
\se \}$ with respect to $\se \in \li [ \f{b}{1 + \vep}, b \ri ]$ is attained at the finite set {\small \[ A_2 =  \li \{ b,  \f{b}{1 + \vep} \ri
\} \bigcup \li \{ \f{k}{n (1 - \vep)} \in \li ( \f{b}{1 + \vep}, b \ri ): k \in \bb{Z} \ri \}.
\]
} \eeL

\bpf

For $\se \in \li ( \f{b}{1 + \vep}, b \ri ]$, it follows from (\ref{ref88c}) of Lemma \ref{lemvipa} that $\Pr \{ | \wt{\bs{\se}}_n - \se | <
\vep \se \mid \se \} = \Pr \{ \wh{\bs{\se}}_n > (1 - \vep) \se \mid \se \}$.  Let $\al < \ba$ be two consecutive elements of $A_2$.  It follows
from  Lemma \ref{plus_rev} that
\[
\Pr \{ | \wt{\bs{\se}}_n - \se | < \vep \se \mid \se \} = \Pr \{ \wh{\bs{\se}}_n > (1 - \vep) \se \mid \se \} = \Pr \{ \wh{\bs{\se}}_n > (1 -
\vep) \se^* \mid \se \}
\]
for $\se \in (\al, \ba)$, where $\se^* = \f{\al + \ba}{2}$.  Now let $\eta \in (0, \f{\ba - \al}{2} )$.  Invoking the assumption that for any
interval $\mscr{I}$, the probability $\Pr \{ Y_n \in \mscr{I} \mid \se \}$ is a unimodal function of $\se \in [a, b]$, we have {\small \bee \Pr
\{ \wh{\bs{\se}}_n
> (1 - \vep) \se \mid \se  \} &  \geq &  \min \li [  \Pr \{ \wh{\bs{\se}}_n > (1 - \vep) (\al + \eta) \mid \al + \eta \}, \; \Pr \{
\wh{\bs{\se}}_n > (1 - \vep) (\ba - \eta) \mid \ba - \eta \} \ri ]\\
& =  & \min \li [  \Pr \{ Y_n \geq g(\al + \eta) \mid \al + \eta \}, \; \Pr \{ Y_n \geq g(\ba - \eta) \mid \ba - \eta \} \ri ] \eee} for any
$\se \in (\al, \ba)$, where the function $g(.)$ is defined by (\ref{ghrev}).  From Lemma \ref{plus_rev}, we know that
\[
g(\al + \eta) = g(\ba - \eta) = g(\al).
\]
Hence, \bee \Pr \{ \wh{\bs{\se}}_n > (1 - \vep) \se \mid \se  \}  \geq   \min \li [  \Pr \{ Y_n \geq g(\al) \mid \al + \eta \}, \; \Pr \{ Y_n
\geq g(\al) \mid \ba - \eta \} \ri ] \eee for any $\se \in (\al, \ba)$.  Recalling the assumption that for any interval $\mscr{I}$, the
probability $\Pr \{ Y_n \in \mscr{I} \mid \se \}$ is a continuous function of $\se \in [a, b]$, we have \bee \Pr \{ | \wt{\bs{\se}}_n - \se | <
\vep \se \mid \se \} & = & \Pr \{ \wh{\bs{\se}}_n > (1 - \vep) \se \mid \se \}\\
 & \geq &  \lim_{\eta \downarrow 0} \min \li [  \Pr \{ Y_n \geq
g(\al) \mid \al + \eta \}, \; \Pr \{ Y_n \geq
g(\al) \mid \ba - \eta \} \ri ]\\
&  = &   \min \li [  \lim_{\eta \downarrow 0} \Pr \{ Y_n \geq g (\al) \mid \al + \eta \}, \; \lim_{\eta \downarrow 0} \Pr \{ Y_n \geq g(\al)
\mid \ba - \eta \} \ri ]\\
&  = &   \min \li [  \Pr \{ Y_n \geq g(\al) \mid \al  \}, \; \Pr \{ Y_n \geq g(\al) \mid \ba  \} \ri ]\\
& \geq &  \min \li [  \Pr \{ Y_n \geq g(\al) \mid \al  \}, \; \Pr \{ Y_n \geq g(\ba) \mid \ba \} \ri ]\\
& = &  \min \li [  \Pr \{ \wh{\bs{\se}}_n > (1 - \vep) \al \mid \al \}, \; \Pr \{ \wh{\bs{\se}}_n > (1 - \vep) \ba \mid \ba \} \ri ]\\
& = &  \min \li [  \Pr \{ | \wt{\bs{\se}}_n - \al | < \vep \al \mid \al \}, \; \Pr \{ | \wt{\bs{\se}}_n - \ba | < \vep \ba \mid \ba \} \ri ]
  \eee for any $\se \in (\al, \ba)$. This immediately leads to the result of the lemma.  \epf

By virtue of (\ref{89933a}) of Lemma \ref{comb33886}, we have the following result.

\beL

\la{lem3388a9}

Assume that $\f{a}{1 - \vep} > \f{b}{1 + \vep}$.  Then, $\Pr \{ | \wt{\bs{\se}}_n - \se | < \vep \se \mid \se \} = 1$ for $\se \in \li ( \f{b}{1
+ \vep}, \f{a}{1 - \vep} \ri )$.  \eeL

By virtue of (\ref{89933b}) of Lemma \ref{comb33886} and a similar argument as that of Lemma \ref{lem3388b96}, we have the following result.

 \beL \la{lem3388b9} Assume that $\f{a}{1 - \vep} > \f{b}{1 +
\vep} \geq a$.   Then, the minimum of $\Pr \{ | \wt{\bs{\se}}_n - \se | < \vep \se \mid \se \}$ with respect to $\se \in \li [a, \f{b}{1 + \vep}
\ri ]$ is attained at the finite set {\small \[ A_3 =  \li \{ a,  \f{b}{1 + \vep} \ri \} \bigcup \li \{ \f{k}{n (1 + \vep)} \in \li (a, \f{b}{1
+ \vep} \ri ): k \in \bb{Z} \ri \}.
\]
} \eeL

By virtue of (\ref{89933c}) of Lemma \ref{comb33886} and a similar argument as that of Lemma \ref{lem3388c96}, we have the following result.

\beL \la{lem3388c9} Assume that $b \geq \f{a}{1 - \vep} > \f{b}{1 + \vep}$.   Then, the minimum of $\Pr \{ | \wt{\bs{\se}}_n - \se | < \vep \se
\mid \se \}$ with respect to $\se \in \li [ \f{a}{1 - \vep}, b \ri ]$ is attained at the finite set {\small \[ A_4 =  \li \{ b,  \f{a}{1 - \vep}
\ri \} \bigcup \li \{ \f{k}{n (1 - \vep)} \in \li ( \f{a}{1 - \vep}, b \ri ): k \in \bb{Z} \ri \}.
\]
} \eeL

We are now in a position to prove the theorem.  We need to consider five cases as follows.

In the case of $\f{a}{1 - \vep} \leq \f{b}{1 + \vep}$, it follows from Lemmas \ref{lem3388a96}--\ref{lem3388c96} that, the minimum of $\Pr \{ |
\wt{\bs{\se}}_n - \se | < \vep \se \mid \se \}$ with respect to $\se \in [a, b]$ is attained at the finite set $( A_0 \cup A_1 \cup A_2 ) \cap
[a, b] = A \cap [a, b]$, where $A$ is the set defined by (\ref{Arev}).

In the case of $b \geq \f{a}{1 - \vep} > \f{b}{1 + \vep} \geq a$, it follows from Lemmas \ref{lem3388a9}--\ref{lem3388c9} that, the minimum of
$\Pr \{ | \wt{\bs{\se}}_n - \se | < \vep \se \mid \se \}$ with respect to $\se \in [a, b]$ is attained at the finite set $( A_3 \cup A_4 ) \cap
[a, b] = A \cap [a, b]$.

In the case of $\f{a}{1 - \vep} > b > \f{b}{1 + \vep} \geq a$, it follows from Lemma \ref{lem3388a9} that $\Pr \{ | \wt{\bs{\se}}_n - \se | <
\vep \se \mid \se \} = 1$ for $\se \in \li ( \f{b}{1 + \vep}, b \ri ]$.  Using this result and Lemma \ref{lem3388b9}, we have that the minimum
of $\Pr \{ | \wt{\bs{\se}}_n - \se | < \vep \se \mid \se \}$ with respect to $\se \in [a, b]$ is attained at the finite set $( A_3 \cup \{b\} )
\cap [a, b] = A \cap [a, b]$.

In the case of $b \geq \f{a}{1 - \vep} > a > \f{b}{1 + \vep}$, it follows from Lemma \ref{lem3388a9} that $\Pr \{ | \wt{\bs{\se}}_n - \se | <
\vep \se \mid \se \} = 1$ for $\se \in \li [ a, \f{a}{1 - \vep} \ri )$.  Using this result and Lemma \ref{lem3388c9}, we have that the minimum
of $\Pr \{ | \wt{\bs{\se}}_n - \se | < \vep \se \mid \se \}$ with respect to $\se \in [a, b]$ is attained at the finite set $( A_4 \cup \{a\} )
\cap [a, b] = A \cap [a, b]$.

In the case of $\f{a}{1 - \vep} > b > a > \f{b}{1 + \vep}$, it follows from Lemma \ref{lem3388a9} that $\Pr \{ | \wt{\bs{\se}}_n - \se | < \vep
\se \mid \se \} = 1$ for $\se \in  [a, b]$.  Thus,  the minimum of $\Pr \{ | \wt{\bs{\se}}_n - \se | < \vep \se \mid \se \}$ with respect to
$\se \in [a, b]$ is attained at the finite set $\{a, b\} \cap [a, b] = A \cap [a, b]$.

The number of elements of the finite set can be calculated by using the property of the ceiling and floor functions. This completes the proof of
the theorem.

\end{document}